\documentclass[preprint,12pt,authoryear]{elsarticle}

\usepackage{amssymb}
\usepackage{amsmath}
\usepackage{amsmath}
\usepackage{subfigure}
\usepackage{url}
\usepackage{multirow} 
\usepackage{booktabs} 
\journal{Elsevier}

\begin{document}

\begin{frontmatter}

%% Title, authors and addresses

%% use the tnoteref command within \title for footnotes;
%% use the tnotetext command for theassociated footnote;
%% use the fnref command within \author or \affiliation for footnotes;
%% use the fntext command for theassociated footnote;
%% use the corref command within \author for corresponding author footnotes;
%% use the cortext command for theassociated footnote;
%% use the ead command for the email address,
%% and the form \ead[url] for the home page:
%% \title{Title\tnoteref{label1}}
%% \tnotetext[label1]{}
%% \author{Name\corref{cor1}\fnref{label2}}
%% \ead{email address}
%% \ead[url]{home page}
%% \fntext[label2]{}
%% \cortext[cor1]{}
%% \affiliation{organization={},
%%            addressline={}, 
%%            city={},
%%            postcode={}, 
%%            state={},
%%            country={}}
%% \fntext[label3]{}

\title{Extracting Interpretable Higher-Order Topological Features across Multiple Scales for Alzheimer's Disease Classification} %% Article title

%% use optional labels to link authors explicitly to addresses:
%% \author[label1,label2]{}
%% \affiliation[label1]{organization={},
%%             addressline={},
%%             city={},
%%             postcode={},
%%             state={},
%%             country={}}
%%
%% \affiliation[label2]{organization={},
%%             addressline={},
%%             city={},
%%             postcode={},
%%             state={},
%%             country={}}

\author[label1]{Dengyi Zhao}
\ead{zhaodengyi@mail.sdu.edu.cn}

\author[label1]{Shanyong Li}
\ead{lishanyong@mail.sdu.edu.cn}

\author[label1]{Yunping Wang}
\ead{wangyunping@mail.sdu.edu.cn}

\author[label1]{Chenfei Wang}
\ead{chenfei_wang0315@163.com}

\author[label2]{Zhiheng Zhou}
\ead{zhouzhiheng@amss.ac.cn}

\author[label2]{Guiying Yan}
\ead{yangy@amt.ac.cn}

\author[label1]{Xingqin Qi\corref{cor1}}
\ead{qixingqin@sdu.edu.cn}

\cortext[cor1]{Corresponding author.}

%\author{Dengyi Zhao$^{1}$ \quad Zhiheng Zhou$^{2}$ \quad Guiying Yan$^{2}$ \quad Dongxiao Yu$^{1}$ \quad  \textbf{Xingqi Qi}$^{1}$\thanks{Corresponding author: \texttt{qixingqin@sdu.edu.cn}} \\
%$^1$Shandong University \quad $^2$University of Chinese Academy of Sciences\\
%\texttt{zhaodengyi@mail.sdu.edu.cn}\\
%\texttt{zhouzhiheng@amss.ac.cn}\\
%\texttt{yangy@amt.ac.cn}\\
%\texttt{dxyu@sdu.edu.cn}\\
%\texttt{qixingqin@sdu.edu.cn}\\

%% Author affiliation
\affiliation[label1]{organization={School of Mathematics and Statistics, Shandong~University}, 
            %addressline={School of Mathematics and Statistics, Shandong University},
            city={Weihai},
            postcode={264209}, 
            state={Shandong},
            country={China}}
\affiliation[label2]{organization={Academy of Mathematics and Systems Science, Chinese Academy of Sciences}, 
            %addressline={School of Mathematics and Statistics, Shandong University},
            city={Beijing},
            postcode={100190}, 
            state={Beijing},
            country={China}}      

%% Abstract

\begin{abstract}
%% Text of abstract
Brain network topology, derived from functional magnetic resonance imaging (fMRI), holds promise for improving Alzheimer's disease (AD) diagnosis. Current methods primarily focus on lower-order topological features, often overlooking the significance of higher-order features such as connected components, cycles, and cavities. These higher-order features are critical for understanding normal brain function and have been increasingly linked to the pathological mechanisms of AD. However, their quantification for diagnosing AD is hindered by their inherent nonlinearity and stochasticity in the brain. This paper introduces a novel framework for diagnosing Alzheimer’s disease that uses persistent homology to extract higher-order topological features from fMRI data. It also introduces four quantitative methods that capture subtle, multiscale geometric variations in functional brain networks associated with AD. Our experimental results demonstrate that this framework significantly outperforms existing methods in AD classification. Extensive ablation studies and interpretability analysis confirm the effectiveness of our framework. Our study also reveals that the number of cycles or cavities significantly decrease in AD patients. The extracted key brain regions derived from cycles and cavities align with domain knowledge in neuroscience literature and provide direct and insightful findings. This study highlights the potential of higher-order topological features for early AD detection and significantly advances the field of brain topology analysis in neurodegenerative disease research.

\end{abstract}
%%Graphical abstract
\begin{graphicalabstract}
\begin{figure}[h]
    \centering
    \includegraphics[scale=0.40]{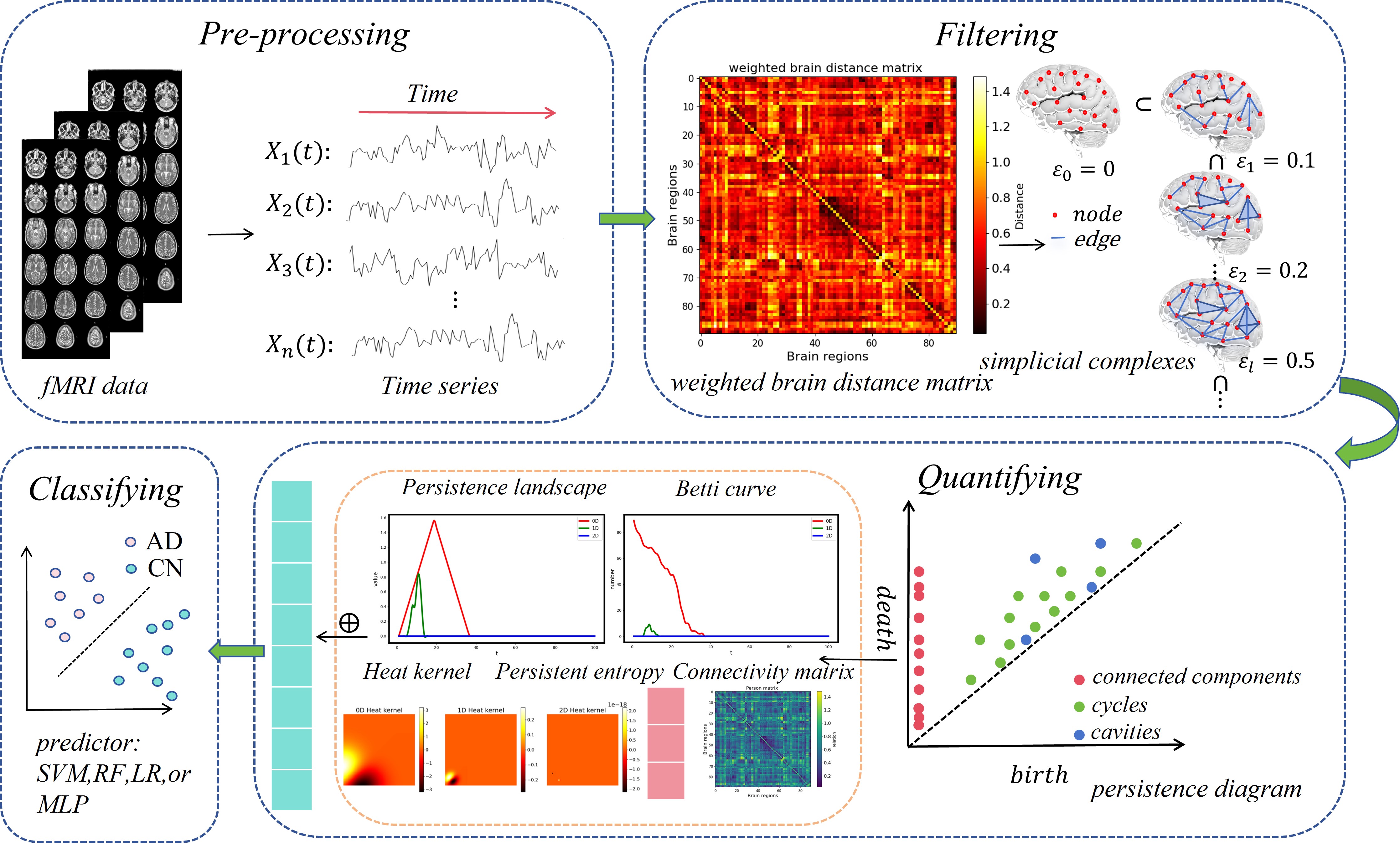}
    %\caption{Enter Caption}
    %\label{fig:enter-label}
\end{figure}
\end{graphicalabstract}

%%Research highlights
\begin{highlights}
\item  Capture multiscale higher-order topological features in fMRI via persistent homology.

\item Introduce four quantitative methods to comprehensively characterize geometric properties of the brain.

\item Reveal that the number of cycles or cavities significantly decrease in AD patients.

\item 
Extracted key brain regions derived from cycles and cavities align with existing  neuroscientific knowledge.
\end{highlights}

%% Keywords
\begin{keyword}
%% keywords here, in the form: keyword \sep keyword
Brain network classification  \sep Higher-order topological features \sep fMRI biomarker \sep Machine Learning

%% PACS codes here, in the form: \PACS code \sep code

%% MSC codes here, in the form: \MSC code \sep code
%% or \MSC[2008] code \sep code (2000 is the default)

\end{keyword}

\end{frontmatter}

%% Add \usepackage{lineno} before \begin{document} and uncomment 
%% following line to enable line numbers
%% \linenumbers

%% main text
%%

%% Use \section commands to start a section
\section{Introduction}
\label{sec1}
%% Labels are used to cross-reference an item using \ref command.

Alzheimer's disease (AD) is a devastating neurodegenerative disorder marked by progressive cognitive decline and memory impairment. The pathogenesis of AD involves a complex interplay of genetic, environmental, psychological, and neurobiological factors, which remains poorly understood and presents significant challenges for early diagnosis and effective treatment \citep{1,2}. Current diagnostic approaches primarily rely on a combination of neuropsychological assessments and clinical observations, which are often inadequate for detecting the disease in its early stages or for making accurate diagnoses \citep{3,4}.

In recent years, functional magnetic resonance imaging (fMRI) has emerged as a promising tool for advancing our understanding and detection of AD \citep{5}. However, fMRI data presents unique challenges due to its high dimensionality, low signal-to-noise ratio, and small sample sizes, making its analysis complex and computationally demanding. A promising line of research focuses on leveraging fMRI data to model the brain as a network, where regions of interest (ROIs) represent nodes and functional connections between them are represented as edges \citep{6}. By analyzing differences in brain network topology between AD patients and cognitively normal (CN) controls, these methods have significantly improved diagnostic accuracy, highlighting the potential of network topology features in unraveling the complexities of AD \citep{7}.

However, current research in this area faces limitations. Traditional approaches often rely on graph theory-based metrics to quantify brain network topology, which are considered lower-order and may not fully capture the intricate neurobiological connectivity patterns \citep{8,9}. Additionally, methods that use graph neural networks (GNNs) usually analyze brain network topology at a fixed scale. This makes them highly dependent on the choice of network construction thresholds \citep{10,11}. This dependence can lead to an oversimplification of the network, potentially introducing bias and causing a significant loss of valuable information. Moreover, it is difficult to associate the extracted latent features with specific organizational structures in the brain. Hypergraph neural network (HGNN) based models, while capable of modeling higher-order topological relationships \citep{12}, are often considered less interpretable. The learned features in HGNNs are embedded in high-dimensional latent spaces, which can also obscure their connection to specific higher-order topological or neurobiological structures characterized by nonlinearity and stochasticity in the brain, similar to GNN-based models \citep{13}. 

\begin{figure}
    \centering
    \includegraphics[scale=0.4]{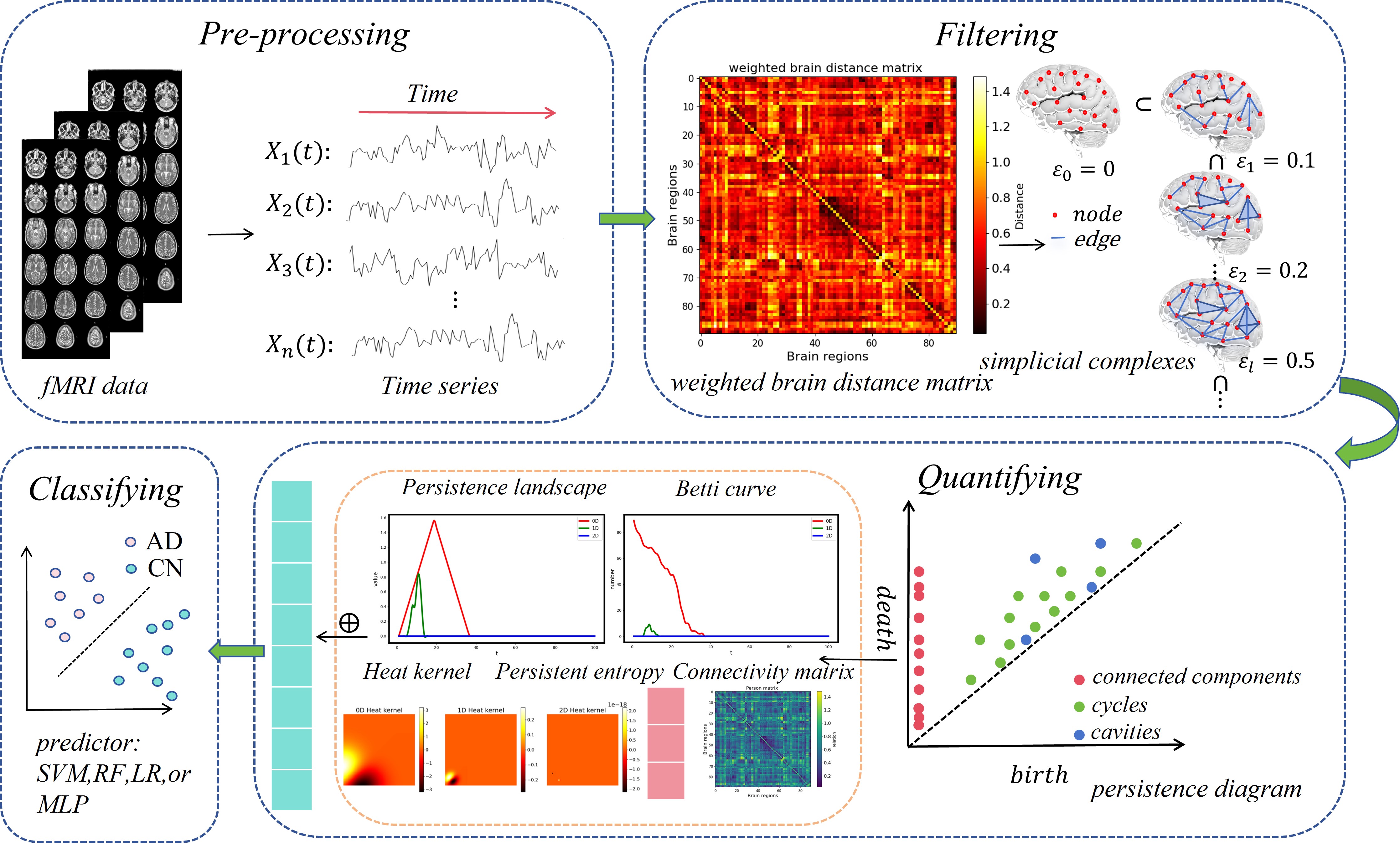}
    \caption{Overview of the proposed framework.
}
    \label{fig1}
\end{figure}

To address these limitations, we propose a novel framework that leverages Persistent Homology \citep{14} to extract both higher-order topological information and lower-order edge features from fMRI data for the diagnosis of brain disorders. This approach captures subtle geometric variations in brain networks across multiple scales and enables a more comprehensive understanding of brain structure and function. Our framework comprises three main stages: filtering, quantification, and classification.
In the filtering stage, fMRI data are encoded as simplicial complexes, a mathematical representation of interconnected brain regions. These are further used to identify higher-order topological features (such as connected components, cycles, and cavities) via Persistent Homology. These features endure across a range of scales, revealing stable geometric characteristics of the brain network at varying levels of detail. They are subsequently summarized in a persistence diagram.
During the quantification stage, four quantitative methods including Persistent landscape, Betti curves, Heat kernels, and Persistent entropy, are employed to extract quantifiable features from the persistence diagram, thereby characterizing the higher-order topological and geometric properties of the brain network from a more comprehensive brain perspective. Additionally, we extract all lower-order edge features from the simplicial complexes with threshold below the maximum distance, capturing lower-order connectivity patterns within the brain.
Finally, in the classification stage, these features are integrated into multiple machine learning models, including Support Vector Machine (SVM), Random Forest (RF), Logistic Regression (LR), and Multilayer Perceptron (MLP), to discriminate between patients with AD and CN. The overall architecture of the proposed framework is illustrated in Figure \ref{fig1}.

In general, this study makes four contributions.
\begin{enumerate}

\item By capturing higher-order topological features in brain networks across multiple scales using Persistent Homology, our method addresses the limitations of fixed-scale network analysis, which can lead to significant loss of valuable topological information.

\item To address the nonlinearity and stochastic characteristics of higher-order features, we employed four quantitative methods (Persistent landscape, Betti curve, Heat kernels, and Persistent Entropy) to capture the higher-order topological and geometric properties of the brain network from a more comprehensive perspective.

\item Our findings show that integrating higher-order topological features with lower-order edge features enables accurate classification of Alzheimer's disease and reveal that the number of cycles or cavities significantly decrease in AD patients.

\item Ablation studies confirm the effectiveness of the framework and interpretability analysis show extracted key brain regions derived from cycles and cavities align with existing neuroscientific knowledge, revealing novel insights for diagnosis of AD.

\end{enumerate}

\section{Related work}
\label{sec2}

\subsection{Traditional Graph Theory Methods}
Traditional graph theory methods have been widely applied to extract topological information from fMRI data in the analysis of brain network structure in Alzheimer’s disease \citep{15}. Previous studies have leveraged graph theoretical approaches to examine the topological characteristics of brain networks in AD patients. For example, Stam et al. \citep{16} investigated abnormal functional network organizations in AD by analyzing small world properties and functional connectivity. Supekar et al. \citep{17} performed a network analysis of intrinsic functional brain connectivity, reporting a loss of small world characteristics in AD patients. Brier et al. \citep{18} found that while characteristic path length remained relatively unchanged in AD compared with CN, clustering coefficient and the number of cliques were significantly reduced. Daianu et al. \citep{19} examined the rich club organization in the AD connectome, identifying structural disruptions in the brain’s core network components. More recently, Yao \citep{20} employed a set of topological metrics, including degree centrality, nodal efficiency, and global efficiency, to quantify network properties, revealing abnormal integration and segregation in the AD group.
Taken together, these studies have primarily focused on lower-order topological features of AD brain networks, which limits their ability to capture higher-order geometric characteristics that may better reflect neurobiologically meaningful connectivity patterns. Such higher-order features, e.g., connected components, cycles, and cavities, play a critical role in supporting normal brain function.

\subsection{Graph Neural Network (GNN)-Based Methods}
In recent years, Graph Neural Networks (GNNs) have emerged as a powerful framework for modeling non-Euclidean topological structures, showing great promise in computer-aided diagnosis of Alzheimer’s disease (AD). Parisot et al. \citep{10} pioneered the application of Graph Convolutional Networks (GCNs) to fMRI data for AD classification. Building on this foundation, Song et al. \citep{21} introduced Dynamic Graph Neural Networks (Dynamic GNNs) to capture the temporal evolution of brain networks in AD. Complementarily, Ma et al. \citep{22} proposed an Attention-Guided Deep Graph Neural Network (AGDGN), which integrates multi-dimensional attention mechanisms with dynamic graph modeling. More recently, Yan et al. \citep{11} developed the Spatio-Temporal Attention Graph Isomorphism Network (STAGIN) to learn latent dynamic representations from fMRI data by incorporating topological information extracted from multiple views (e.g., Pearson correlation, mutual information).
Despite these advances, existing GNN-based approaches often struggle to explicitly disentangle lower-order from higher-order topological information during learning, thereby limiting their ability to capture the intricate non-linear relationships embedded in brain networks. Furthermore, the performance of graph neural models is highly dependent on the quality of network construction and the informativeness of node features; oversimplification of the network or insufficiently descriptive features can substantially reduce diagnostic accuracy. Finally, the limited interpretability of most graph-based frameworks obscures the underlying neurobiological mechanisms driving individual predictions.

\subsection{Higher-Order Topological Information Extraction Methods}
The analysis of brain networks has increasingly emphasized higher-order features that extend beyond simple pairwise connections to reveal more intricate organizational patterns \citep{23}. Recent advancements have introduced several methods for extracting such features. Wang et al. \citep{24} proposed an evolving hypergraph convolutional network for dynamic hyper-brain networks, incorporating an attention mechanism to enhance representation learning and facilitate AD diagnosis. Similarly, Yang et al. \citep{25} constructed higher-order functional connectivity networks using hypergraphs, aiming to capture genuine interaction patterns among brain regions. However, hypergraph neural network approaches often embed features in high-dimensional spaces, making them difficult to interpret in biological or topological terms. Moreover, these methods typically require predefined hypergraph structures, which limits their ability to adapt to the dynamic and scale-dependent changes of brain networks.

To overcome these limitations, Lee et al. \citep{26} pioneered the application of Persistent Homology (PH), an algebraic topology tool, to brain network analysis. They introduced a framework for extracting higher-order topological features from electroencephalogram (EEG) data, addressing the arbitrary thresholding issue inherent in traditional graph-theoretical methods. Nevertheless, the complexity and noise of fMRI data pose substantial challenges for directly applying this method to Alzheimer’s disease. Cassidy et al. \citep{27} leveraged PH to examine changes in connected components of brain networks, focusing on the stability of functional connectivity across spatial and temporal scales. Despite its utility, their approach did not consider higher-dimensional structures, such as cycles, which could provide additional insights. Sizemore et al. \citep{28} extended this line of work by identifying higher-order topological features, including cliques and cavities, in brain networks using algebraic topology, thereby moving beyond pairwise interactions to uncover more complex relationships. More recently, Bian et al. \citep{29} proposed an adversarial training-based Persistent Homology Graph Convolutional Network (ATPGCN) to capture disease-specific connectome patterns and classify brain disorders. This approach integrates loop and component features derived from algebraic topological analysis with readout features obtained from the global pooling layer of GCNs, enabling collaborative learning of individual-level representations.

However, existing studies have not comprehensively quantified higher-order topological features extracted via PH, overlooking valuable information due to their nonlinear and stochastic nature, which complicates quantification. To address this challenge, the Persistence Landscape method \citep{30} was introduced, which transforms complex nonlinear structures into piecewise linear functions in a Banach space, thereby simplifying the quantification of high-dimensional PH features. Other methods have also been developed for quantifying such features from a more comprehensive perspective, including Betti curves \citep{32}, Heat kernel \citep{33}, and Persistent entropy \citep{31}.

%% Use \subsubsection, \paragraph, \subparagraph commands to 
%% start 3rd, 4th and 5th level sections.
%% Refer following link for more details.
%% https://en.wikibooks.org/wiki/LaTeX/Document_Structure#Sectioning_commands

\section{Method}
\subsection{Preliminaries}
\label{3.1}
\subsubsection{Notations}

We employ a standard brain atlas to partition the brain into \(N\) regions of interest (ROIs), with each node in the brain network corresponding to a distinct region. Let \(X_i(t)\) denote the observed blood-oxygen-level–dependent (BOLD) signal from brain region \(i \in V\) at time \(t \in \{1,\dots,T\}\), where \(V = \{1,\dots,N\}\) represents the set of all brain regions. We define \(W\) as the weighted brain distance matrix, where each element \(W_{ij} = {d_{X}}(X_i, X_j)\) quantifies the dependence between regions \(i\) and \(j\). Furthermore, we introduce a network filter, \(\mathcal{B}(X, \varepsilon)\), which produces a thresholded network in which nodes correspond to brain regions, and an edge is established between nodes \(i\) and \(j\) if their distance \(d_X(X_i, X_j)\) falls below a predefined threshold \(\varepsilon\).

\subsubsection{Problem Formulation}

The objective of this study is to formulate a classification problem analogous to graph classification. Specifically, we aim to determine whether a subject has Alzheimer’s disease based on a collection of brain network graphs,  
\[
\mathcal{B} = \bigl\{B_0(X,\varepsilon_0),\, B_1(X,\varepsilon_1),\, B_2(X,\varepsilon_2),\, \dots,\, B_n(X,\varepsilon_n)\bigr\}, \tag{1}
\]
which are constructed at multiple thresholds \(\varepsilon\) from fMRI time series data. The central objective of the proposed model is to extract higher-order topological features from these brain networks and subsequently utilize them for classification.

\subsection{Construction of Weighted Brain Distance Matrix}

The preprocessed fMRI time series signals are used as input, with Pearson correlation selected as the similarity measure between brain regions. Based on this measure, we construct a connectivity matrix from the fMRI data. To account for the characteristics of the time series signals, we derive a weighted brain distance matrix \(W\). Specifically, we compute a pairwise distance matrix between brain regions, where the dissimilarity between signals at regions \(i\) and \(j\) is defined as

\[
W_{ij} = d_X(X_i, X_j) = 1 - \frac{\operatorname{cov}(X_i, X_j)}{\sigma_{X_i}\sigma_{X_j}}, \tag{2}
\]
where \(\operatorname{cov}(X_i, X_j)\) denotes the covariance between regions \(i\) and \(j\), and \(\sigma_{X_i}\) and \(\sigma_{X_j}\) are the standard deviations of the corresponding signals. The weighted brain distance matrix is then employed to construct simplicial complexes, which provide a higher-dimensional generalization of graphs and encode the relationships among brain regions. By transforming the brain network into a sequence of simplicial complexes, we apply persistent homology to extract informative topological features.

\subsection{Filtration of Brain Networks}

In brain network analysis, the choice of a threshold parameter can substantially influence the results, often leading to oversimplification of the network, the introduction of bias, and the loss of valuable information. To address this challenge, we employ the Vietoris-Rips filtration method. This approach constructs a sequence of simplicial complexes that progressively approximate the weighted brain connectivity network with increasing levels of detail. Formally, we define a filtration of simplicial complexes as  
\[
B_0(X, \varepsilon_0) \subset B_1(X, \varepsilon_1) \subset \ldots \subset B_l(X, \varepsilon_l) \subset \ldots \subset B_n(X, \varepsilon_n), \tag{3}
\]  
where each complex \(B_l(X, \varepsilon_l)\) is contained within the subsequent one. The filtration is governed by a scanning parameter \(\varepsilon_l\), which increases from 0 to 2. At each step \(l\), we construct a weighted simplicial complex \(B_l(X, \varepsilon_l)\), where an edge is included if the distance is smaller than \(\varepsilon_l\).

\subsection{Higher-Order Feature Extraction of Brain Networks}
\label{3.2}

During the Rips filtration process, a $k$-dimensional hole (a higher-order topological feature) is born at a filtration value $b_i$ and dies at $d_i$, with its persistence given by the interval $(b_i, d_i)$. These features are higher-order topological invariants that characterize the shape of a space: $0$-dimensional holes correspond to connected components, $1$-dimensional holes to one-dimensional cycles, and $2$-dimensional holes to two-dimensional cavities. The persistence of these higher-order features is summarized in a persistence diagram (PD), defined as the multiset of birth--death pairs:
\[
{PD} = \{(b_1, d_1), (b_2, d_2), \ldots, (b_m, d_m)\}.  \tag{4}
\]
This diagram is a scatter plot with birth values on the x-axis and death values on the y-axis. Since $b_i < d_i$ for all non-trivial features, all points lie above the diagonal line $y = x$. The persistence diagram provides a stable summary of the network's multiscale topology, offering significant insight into its underlying structure and function.

\subsection{Quantification of Higher-Order Topological Features}
\label{quantification of Persistent Topological Features}
Building on previous steps, we generated persistence diagrams that reflect the persistence of specific topological features within an individual's brain network. These diagrams were then transformed into feature vectors suitable for input into machine learning models. To achieve this, we employed four distinct methods: Persistence landscape \citep{30}, Betti curve \citep{32}, Heat kernel \citep{33}, and Persistent entropy \citep{31}.

\textbf{1) Persistence landscape:} The persistence landscape maps points from a persistence diagram to a function space, preserving the core information of higher-order topological features (birth, death, and persistence) and making it usable by traditional statistical and machine learning tools. Specifically, for each point $(b_i, d_i)$ in the $k$-dimensional persistence diagram $PD_k = \{(b_1, d_1),\ldots ,(b_i, d_i)\}_{k}$ (where $b_i$ and $d_i$ are the birth and death times of the $i$-th $k$-dimensional topological feature, respectively), we define a landscape function $f_i(t)$:
\[
f_i(t) =
\begin{cases}
t - b_i & \text{if } t \in \left[ b_i, \frac{b_i + d_i}{2} \right], \\
d_i - t & \text{if } t \in \left[ \frac{b_i + d_i}{2}, d_i \right], \\
0 & \text{otherwise.}
\end{cases}
\tag{5}
\]
Here, $t$ is a point in the domain of the function, which can be interpreted as a continuous variable such as filtration time. This piecewise linear function resembles a mountain peak.

\textbf{2) Betti curve:}
Betti numbers are topological invariants that quantify the number of holes at different dimensions within a space. The Betti curve is constructed by calculating the Betti number at various filtration values, providing a dynamic view of the persistence of topological features within the brain network. For each scale $t$, the Betti number is calculated as:
\[
\beta_k(t) = \#\{(b_i, d_i) \in PD_k \mid b_i \leq t < d_i\},\tag{6}
\]
where $\#$ denotes the cardinality of the set.

\textbf{3) Heat kernel:}
The Heat kernel method captures fine-grained feature information from persistence diagrams through Gaussian kernel convolution. This process enables a precise characterization of the topological feature distributions and their evolution. The Heat kernel method generates feature vectors by treating the persistence diagram as a heat diffusion process. A Heat kernel function is constructed and then transformed into feature vectors. Specifically, for the $k$-dimensional persistence diagram of topological features, we define a Heat kernel function for each point $(b_i, d_i)$ in the diagram:
\[
K_i(x, y; \sigma) = \frac{1}{2\pi\sigma^2} \exp\left(-\frac{(x-b_i)^2 + (y-d_i)^2}{2\sigma^2}\right), \tag{7}
\]
where $\sigma$ is the Gaussian bandwidth controlling the diffusion scale.

\textbf{4) Persistent Entropy:}
Persistence entropy is an information-theoretic approach that quantifies the uniformity and complexity of topological feature distributions in persistence diagrams. For the $k$-dimensional persistence diagram of topological features, we first calculate the persistence duration $l_i = d_i - b_i$ for each point $(b_i, d_i)$, and then normalize these durations to form a probability distribution:
\[
p_i = \frac{l_i}{\sum_{j=1}^{n} l_j} .\tag{8}
\]
The persistence entropy feature vector for $k$-dimensional topological features is defined as:
\[
pe_k = -\sum_{i=1}^{n} p_i \log (p_i). \tag{9}
\]
Higher persistence entropy indicates a more uniform feature distribution across scales and greater structural complexity, while lower values suggest a concentration of features at specific scales. In Alzheimer's disease research, changes in persistence entropy serve as an indicator of disrupted brain network topology.

\subsection{Lower-Order Feature Extraction of Brain Networks}

It should be noted that all the aforementioned features represent higher-order topological characteristics of brain networks. However, lower-order information plays an equally important role in subsequent feature classification. 
Specifically, this features correspond to extracting all edge features with distances below the maximum threshold from the brain network graphs $\{B_n(X,\varepsilon_n)\}$.

\subsection{Classification by Machine Learning}
The final fused feature vector $\mathbf{v}$ integrates multi-scale topological characteristics by concatenating both higher-order and lower-order features:
\begin{equation}
\mathbf{v} = \mathbf{lc} \oplus \mathbf{bc} \oplus \mathbf{heat} \oplus \mathbf{pe} \oplus \mathbf{c}
\tag{10}
\end{equation}
where $\mathbf{lc}$, $\mathbf{bc}$, $\mathbf{heat}$, and $\mathbf{pe}$ represent the higher-order persistent homology features (Persistence landscape, Betti curve, Heat kernel, and Persistent entropy, respectively), while $\mathbf{c}$ captures lower-order connectivity patterns.

This vector $\mathbf{v}$ can be used as input for various machine learning classification models. Traditional statistical methods struggle to analyze such complex data, whereas machine learning algorithms can automatically identify the relationships between these features and the presence of AD, leading to more accurate and generalizable classification results. In this study, we employed several machine learning models, including Support Vector Machine (SVM) \citep{34}, Random Forest (RF) \citep{35}, Logistic Regression (LR) \citep{36}, and Multi-Layer Perceptron (MLP) \citep{34}, to classify AD based on the extracted topological features. These models were selected for their proven effectiveness in handling complex, high-dimensional data and their ability to provide reliable predictions in clinical settings.

\section{Experiments}
\subsection{Experimental Settings}
\subsubsection{Datasets and Preprocessing}
We employed 232 functional magnetic resonance imaging (fMRI) datasets from the ADNI1 database (https://adni.loni.usc.edu/) \citep{39}, consisting of 116 Alzheimer's disease (AD) cases and 116 cognitively normal (CN) subjects. The fMRI data were preprocessed using the SPM pipeline \citep{38}. Specifically, the T1-weighted (T1w) images were used for brain masking, image alignment, and BOLD signal normalization. The preprocessing steps included motion correction, realignment, field unwarping, normalization, bias field correction, and brain extraction. Subsequently, the AAL atlas \citep{40} was applied to parcellate the brain into 90 regions of interest (ROIs). The BOLD time series of each ROI was obtained by averaging the signals across all voxels within the region. For each subject, a weighted connectivity matrix was constructed, where each element represents the Pearson correlation between the average BOLD signals of paired ROIs.

\subsubsection{Methods for Comparison}
To evaluate the effectiveness of our approach, we compared it against ten baseline methods, grouped into three categories: 

1) \textbf{Traditional graph
theory-based models}: betweenness centrality (BC), local efficiency (LE), clustering coefficient (CCO), and the upper-triangular elements of the connectivity matrix (UTE) \citep{41}. These features were used as inputs to a Multi-Layer Perceptron (MLP) \citep{34} for classification.  

2) \textbf{Graph neural network (GNN)-based models}: GCN \citep{42}, GAT \citep{44}, GraphSAGE \citep{43}, and Graph Transformer \citep{45}.  

3) \textbf{Higher-order models}: HGCN \citep{46} and HGAT \citep{47}.

\subsubsection{Implementation Details}
To ensure robust evaluation, we employed a 5-fold cross-validation strategy for model training and testing, and reported the average performance across folds as the final result. Accuracy, Precision, Recall, and F1-score were used as evaluation metrics.  

Higher-order topological features were extracted using four functions from the giotto-tda library \citep{48}: \texttt{Persistence landscape} (discretized with $n_{\text{bins}}=100$ and $n_{\text{layers}}=1$ or $2$, yielding $3 \times 100$ and $6 \times 100$ vectors), \texttt{Betti curve} ($n_{\text{bins}}=100$, producing $3 \times 100$ vectors), \texttt{Heat kernel} (with $\sigma=1.2$ and $1.4$, generating two $3 \times 10^2$ matrices), and \texttt{Persistent Entropy} ($3 \times 1$ vectors). In addition, unfolding the upper-triangular part of each connectivity matrix produced a vector of dimension 4005. All features were concatenated into a $5808 \times 1$ feature vector for machine learning analysis.  

We implemented four classifiers: a Support Vector Machine (SVM) with an RBF kernel ($\gamma=0.001$, $C=0.8$); a Random Forest (RF) with 300 trees and a maximum depth of 20; a Logistic Regression (LR) model with regularization strength $C=3.6$, using the \texttt{lbfgs} solver and a maximum of 500 iterations; and a Multilayer Perceptron (MLP) with three hidden layers, L2 regularization, and ReLU activation. All other hyperparameters were kept at their default values to ensure fair comparison. This configuration balanced model complexity and regularization while effectively leveraging the high-dimensional feature space.

\begin{table}[ht]
  \centering
  \caption{Performance comparison with four categories of baselines on the ADNI dataset. The best results are marked in bold and the second-best results are underlined.}
  \scalebox{0.84}{
  \begin{tabular}{@{}cccccc@{}}
    \toprule
    \multirow{2}{*}{Type}  & \multirow{2}{*}{Method} & \multicolumn{4}{c}{ADNI}  \\
    \cmidrule(lr){3-6} 
    & & Accuracy & Precision & Recall & F1-score \\
    \midrule
    \multirow{5}{*}{\parbox{3cm}{\centering Traditional Graph Theory-Based Models}}
    & BC+MLP  & 0.625 & 0.615 & 0.761 & 0.680 \\
    & LE+MLP  & 0.750 & 0.777 & 0.700 & 0.736 \\
    & CCO+MLP  & 0.725 & 0.695 & 0.800 & 0.744 \\
    & UTE+MLP  & 0.775 & 0.739 & 0.850 & 0.790 \\
    & LE+UTE+MLP & 0.800 & 0.772 & 0.850 & 0.809 \\
    \midrule
    \multirow{4}{*}{\parbox{3cm}{\centering GNN-Based Models}}
   & GCN  & 0.793 & 0.763 & 0.856 & 0.803 \\
    & GAT  & 0.829 & 0.761 & 0.842 & 0.800 \\
    & GraphSAGE  & 0.851 & 0.857 & 0.889 & 0.872 \\
    & Graph Transformer  & 0.829 & \underline{0.947} & 0.720 & 0.818 \\
    \midrule
    \multirow{2}{*}{\parbox{3cm}{\centering HGNN-Based Models}}
    & HGCN  & 0.815 & 0.857 & 0.907 & 0.800 \\
    & HGAT  & 0.850 & 0.825 & 0.880 & 0.831 \\
   \midrule
    \multirow{4}{*}{\parbox{3cm}{\centering our framework}} 
   & Higher-order+SVM & \textbf{0.936} & \textbf{0.957} & \underline{0.917} & \textbf{0.936} \\
    & Higher-order+RF & 0.851 & 0.833 & 0.870 & 0.851 \\
    & Higher-order+LR & \underline{0.915} & 0.880 & \textbf{0.957} & \underline{0.917} \\
    & Higher-order+MLP  & \textbf{0.936} & 0.917 & \textbf{0.957} & \textbf{0.936} \\
    \midrule

  \end{tabular}
  \label{tab:performance_comparison on adni} % Add label for referencing
  }
\end{table}

\subsection{Model Comparison}
In this section, we compare the performance of our proposed model with baseline methods for Alzheimer's disease diagnosis. The detailed results are reported in Table \ref{tab:performance_comparison on adni}, where the best outcomes for each metric are shown in bold and the second-best results are underlined. Our higher-order topological feature extraction framework consistently achieves superior performance across all four evaluation metrics. Specifically, compared with traditional graph theory-based feature extraction methods, our model improves accuracy by 20--30\%. Relative to Graph Neural Network (GNN)-based models, our framework yields an accuracy gain of approximately 13\%, and compared with higher-order topological information extraction models, it achieves an additional improvement of about 10\%. Moreover, models that integrate higher-order topological features with machine learning classifiers outperform graph theory-based, GNN-based, and HGNN-based approaches in terms of precision, recall, and F1-score.  

Traditional graph theory-based methods consistently yield lower performance across all evaluation metrics, underscoring the effectiveness of modeling brain networks with higher-order topological features. Although GNN-based models can capture topological structures, their performance remains below that of our framework, suggesting that they primarily extract lower-order information and fail to adequately represent complex, higher-order nonlinear relationships between brain regions. Higher-order topological information extraction models achieve better results than GNN-based approaches but still underperform our method. This limitation arises because the high-level information they capture is difficult to associate with specific organizational structures in the brain, thereby reducing classification accuracy.

Among the combined models, Higher-order+SVM achieves the highest accuracy and precision, whereas Higher-order+MLP provides the best recall and F1-score. These results highlight the importance of integrating higher-order topological features with advanced machine learning models to achieve state-of-the-art performance in AD classification.

\subsection{Model Analysis}
\subsubsection{Ablation Study}
In this section, we evaluate the effectiveness of the design choices for different components of the model. Specifically, we examine their impact on experimental outcomes, including the representation of lower-order features (the upper-triangular similarity matrix represented as vector \textit{c}), the selection of quantification methods for higher-order persistent topological features—such as persistence landscape, Betti curve, Heat kernel, and persistenct entropy—as well as the influence of feature dimensionality.  

\begin{figure}[ht]
    \centering
    \includegraphics[scale=0.4]{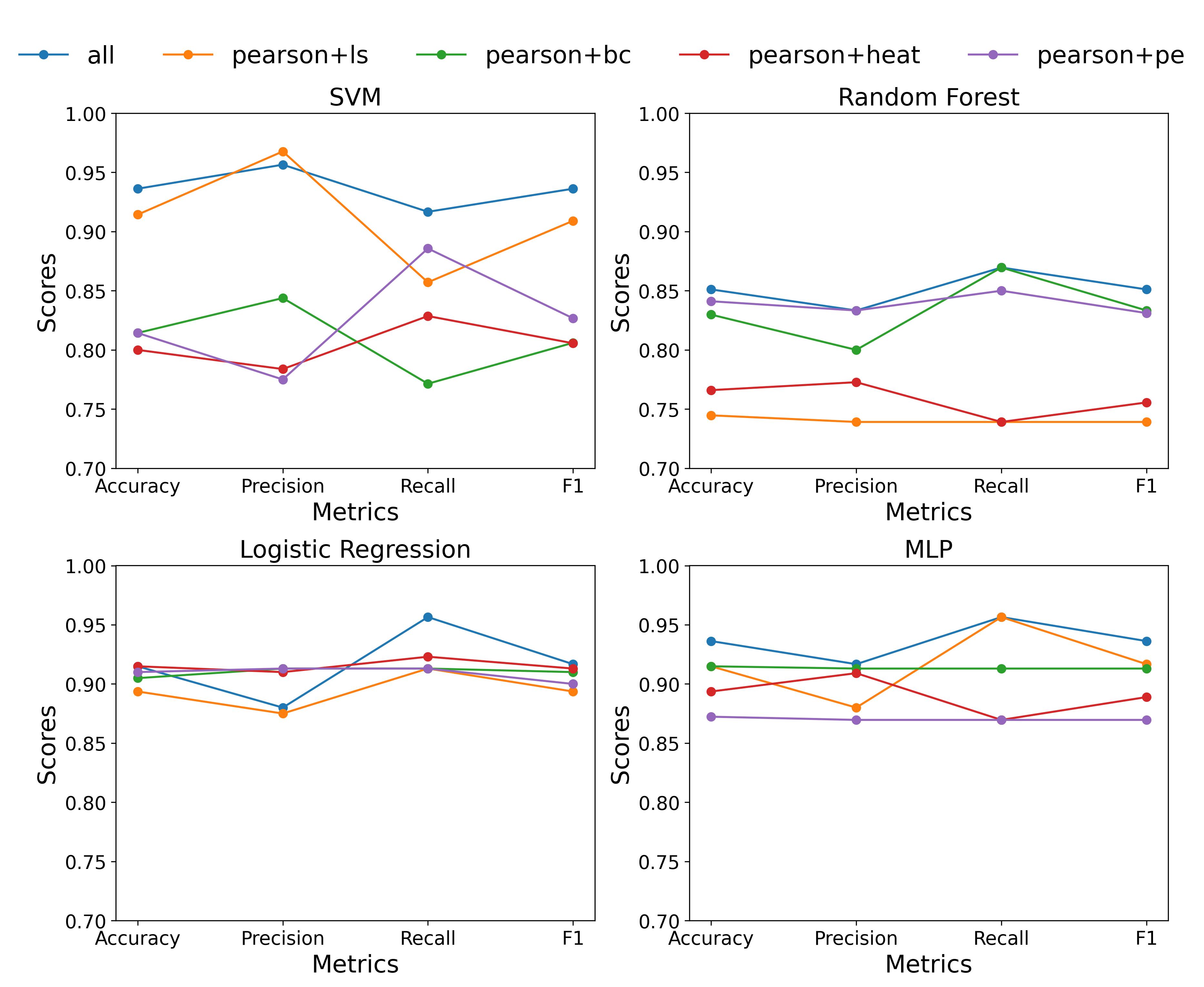}
    \caption{The results of ablation study on vectorized higher-order persistent topological features.}
    \label{fig2}
\end{figure}

\textbf{1) Quantification of Persistent Topological Features: }Four distinct vectorization methods for persistence diagrams were employed to extract higher-order topological features. To assess the individual contribution of each method, we conducted comparative experiments using four machine learning algorithms: SVM, RF, LR, and MLP. As shown in Figure \ref{fig2}, the framework achieves the best overall performance when all four vectorization methods are combined. Notably, while the integrated approach provides the highest classification accuracy across most evaluation metrics, certain individual methods achieve competitive results in specific scenarios. For instance, the persistence landscape method alone yields relatively high precision in SVM models, and the combination of the Heat kernel with Betti curves demonstrates comparable performance in LR models. However, these individual methods underperform on other metrics, as they primarily capture specific aspects of topological information while lacking comprehensive coverage of higher-order brain network properties. These findings suggest that all four vectorization approaches—Persistence landscape, Betti curve, Heat kernel, and Persistent entropy—provide complementary and non-redundant information for characterizing the complex topological features of brain networks in Alzheimer's disease classification.

\begin{figure}[ht]
    \centering
    \includegraphics[scale=0.4]{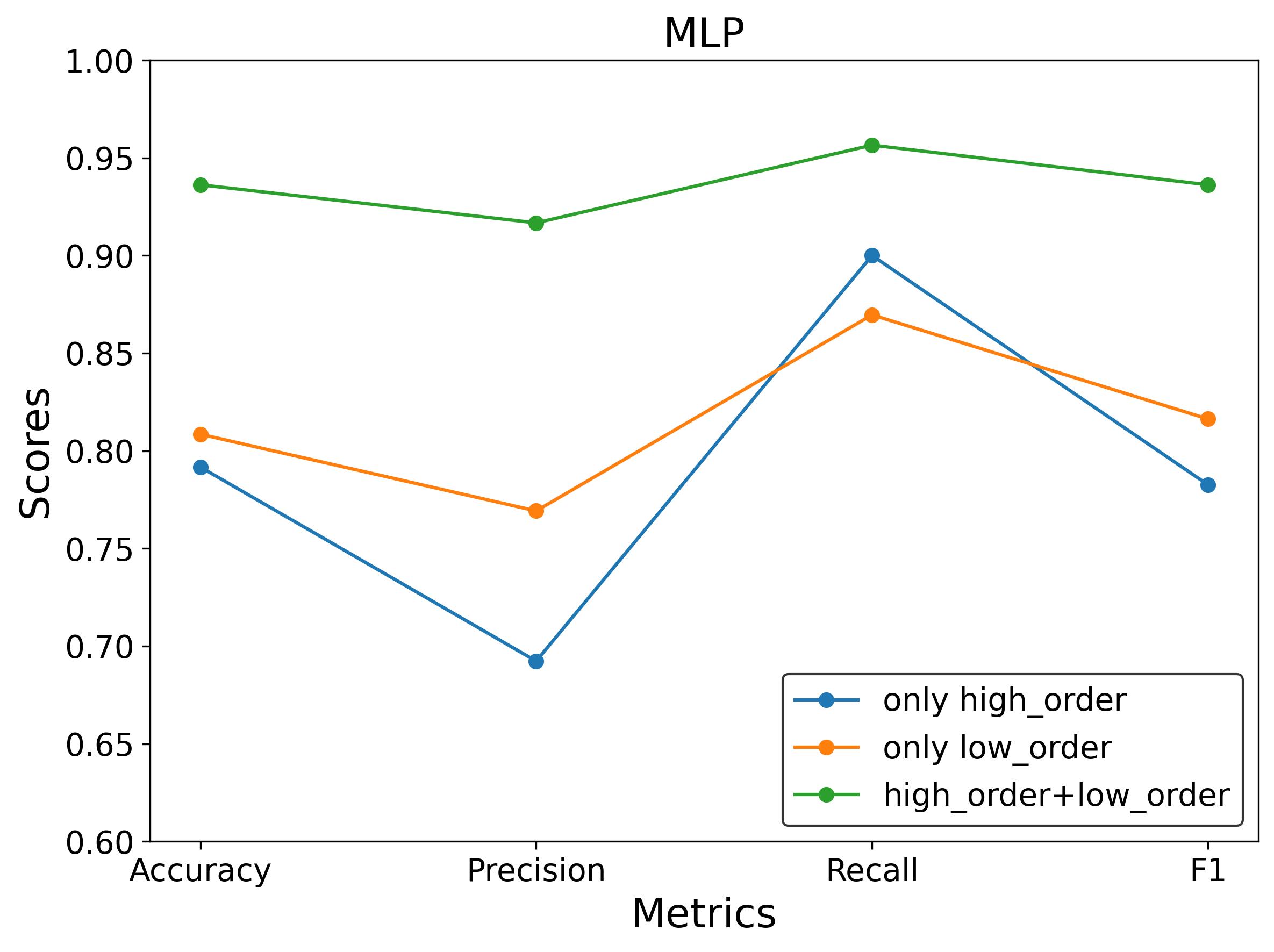}
    \caption{The results of ablation study on higher-order and lower order topological features.}
    \label{fig3}
\end{figure}

\begin{figure}[ht]
    \centering
    \includegraphics[scale=0.4]{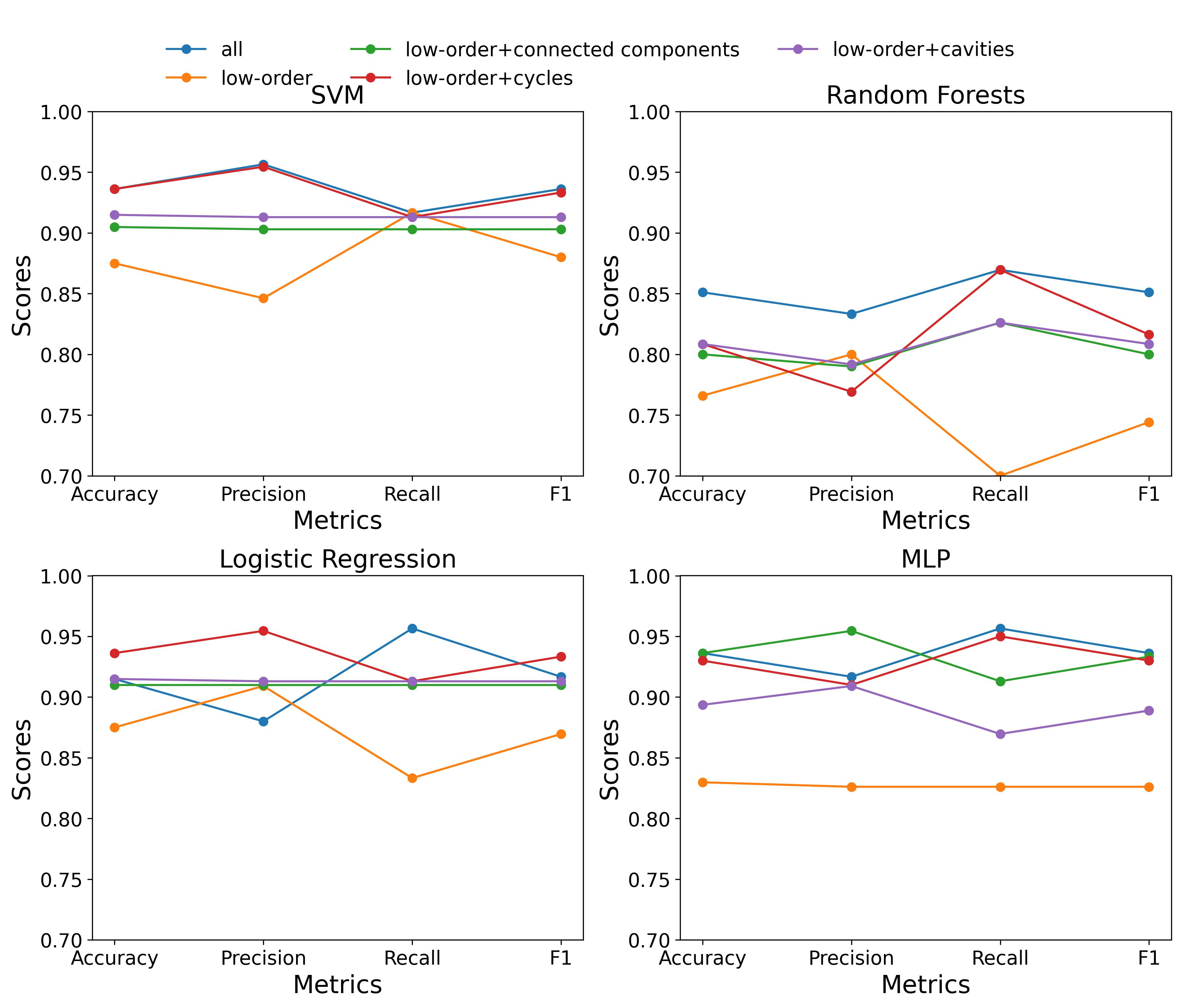}
    \caption{The results of ablation study on selection of higher-order topological features}
    \label{fig4}
\end{figure}

\textbf{2) Higher-Order and Lower-Order Topological Features: } 
We conducted experiments to examine the respective contributions of higher-order persistent topological features and lower-order topological features, thereby clarifying their relative importance. The results, shown in Figure \ref{fig3}, demonstrate that the model achieves the best overall performance when both feature types are fused. In contrast, using only higher-order or only lower-order features for classification yields inferior results across all four evaluation metrics. This finding highlights the limitations of analyzing Alzheimer's patients' brain network topology from a single lower-order perspective, as such an approach fails to capture the full complexity of the underlying structures.  

\textbf{3) Selection of Higher-Order Topological Features: } 
To evaluate the contributions of different higher-order topological features—connected components, cycles, and cavities—we conducted ablation experiments by selectively excluding each feature during persistent homology analysis. In particular, we assessed performance when lower-order features were fused with higher-order features from individual dimensions. As shown in Figure \ref{fig4}, the model achieves the best performance when lower-order features are combined with higher-order features from all dimensions. These results indicate that connected components, cycles, and cavities provide complementary information for characterizing brain network topology. Together, they comprehensively describe higher-order structural properties from multiple perspectives. Notably, even when only connected components or cavities are fused with lower-order features, the classification performance still exceeds that of traditional graph-theoretical metrics. This finding underscores the essential role of higher-order persistent topological information in the diagnosis of Alzheimer's disease.

\subsubsection{Hyperparameter Analysis}
This section analyzes the sensitivity of our best-performing model, Higher-order+MLP, to the key hyperparameters of three higher-order topological feature vectorization methods: Persistent landscape, Betti curve, and Heat kernel. Persistent Entropy is excluded, as it does not involve critical hyperparameters.

\textbf{1) Influence of hyperparameter $m_1$ in Persistent landscape: }
The filtration parameter $m_1$, which determines the number of sampling points across homology dimensions, plays a crucial role in quantifying the persistent landscape. Specifically, $m_1$ controls the resolution at which dynamic changes in topological features are captured. A small $m_1$ may lead to insufficient sampling and loss of multiscale information, whereas a large $m_1$ increases computational cost and can introduce noise, thereby degrading performance. 
To assess its impact, we conducted experiments with $m_1$ ranging from 50 to 200 (see Table \ref{tab:landscape}). The results show that setting $m_1=100$ yields optimal classification accuracy, striking an effective balance between computational efficiency and accurate characterization of scale-dependent topological dynamics. Importantly, our model demonstrates low sensitivity to variations in $m_1$, ensuring stable and robust performance across different parameter configurations.

\begin{table}[htbp]
\centering
\caption{Experimental results of the hyperparameter analysis for Persistent landscape.}
\label{tab:landscape}
\scalebox{0.85}{
\begin{tabular}{ccccc}
\toprule
\(m_1\) & Accuracy & Precision & Recall & F1-score \\
\midrule
50    & 0.7917   & 0.6923   & 0.9000   & 0.7826   \\
100   & 0.8333   & 0.8000   & 0.8000   & 0.8000   \\
150   & 0.7500   & 0.6667   & 0.8000   & 0.7273   \\
200   & 0.7500   & 0.6429   & 0.9000   & 0.7500   \\
\bottomrule
\end{tabular}
}
\end{table}

\textbf{2) Influence of hyperparameter $m_2$ in Betti curve:} 
In Betti curve quantification within persistent homology analysis, the number of sampled filtration parameter values ($m_2$) per homology dimension is a key factor in determining the resolution of topological feature characterization. Specifically, $m_2$ governs the ability of Betti curves to capture multiscale topological variations in brain networks. 
A small $m_2$ may lead to under-sampling, resulting in an incomplete representation of the dynamic evolution of topological features and reduced discriminatory power between Alzheimer's disease patients and healthy controls. Conversely, a large $m_2$ increases computational complexity and may introduce excessive detail, which can compromise model stability and generalization. 
To evaluate this effect, we tested $m_2$ values ranging from 50 to 200 (see Table \ref{tab:betti}). The results show that setting $m_2=100$ achieves peak classification performance, enabling Betti curves to effectively track quantitative changes in multidimensional topological features, such as connected components, cycles, and cavities, while maintaining manageable computational cost.

\begin{table}[htbp]
\centering
\caption{Experimental results of the hyperparameter analysis for Betti curve.}
\label{tab:betti}
\scalebox{0.85}{
\begin{tabular}{ccccc}
\toprule
\(m_2\) & Accuracy & Precision & Recall & F1-score \\
\midrule
50    & 0.7083   & 0.6154   & 0.8000   & 0.6957   \\
100   & 0.8333   & 0.8000   & 0.8000   & 0.8000   \\
150   & 0.7083   & 0.6154   & 0.8000   & 0.6957   \\
200   & 0.7083   & 0.6000   & 0.9000   & 0.7200   \\
\bottomrule
\end{tabular}
}
\end{table}

\textbf{3) Influence of hyperparameter $\sigma$ in Heat kernel:} 
The standard deviation parameter ($\sigma$) of the Gaussian kernel is a critical factor in Heat kernel methods for persistent homology analysis, as it controls the resolution of topological feature distributions in persistence diagrams. A small $\sigma$ results in a narrow bandwidth that captures only localized features while neglecting broader structural patterns, leading to incomplete representations of network dynamics. Conversely, a large $\sigma$ produces excessive smoothing that obscures important topological details and reduces discriminative power for Alzheimer's disease classification. 
We systematically evaluated $\sigma$ values from 1.0 to 2.0 (see Table \ref{tab:heatkernel}). The results indicate that $\sigma=1.2$ provides the optimal balance, preserving essential topological features while applying sufficient smoothing. This configuration achieves superior performance in capturing the multiscale evolution of topological characteristics in brain networks, while maintaining computational efficiency and ensuring accurate representation of disease-relevant patterns.

\begin{table}[htbp]
\centering
\caption{Experimental results of the hyperparameter analysis for Heat kernel.}
\label{tab:heatkernel}
\scalebox{0.85}{
\begin{tabular}{ccccc}
\toprule
\(\sigma\) & Accuracy & Precision & Recall & F1-score \\
\midrule
1.0   & 0.7500   & 0.7000   & 0.7000   & 0.7000   \\
1.2   & 0.8333   & 0.7500   & 0.9000   & 0.8182   \\
1.4   & 0.8333   & 0.8000   & 0.8000   & 0.8000   \\
1.6   & 0.7917   & 0.7273   & 0.8000   & 0.7619   \\
1.8   & 0.7917   & 0.7778   & 0.7000   & 0.7368   \\
2.0   & 0.7500   & 0.6429   & 0.9000   & 0.7500   \\
\bottomrule
\end{tabular}
}
\end{table}

\subsection{Interpretability Analysis}
To further investigate the mechanisms underlying the improved classification performance of our model, we designed three interpretability experiments: (1) assessing group differences in higher-order topological features between healthy controls and patients; (2) examining whether the latent features extracted by the model exhibit clear separation between the two groups; and (3) identifying potential biomarkers that align with findings reported in the literature.

\subsubsection{Interpretability of Higher-Order Topological Features}
\begin{figure}
    \centering
    \includegraphics[scale=0.3]{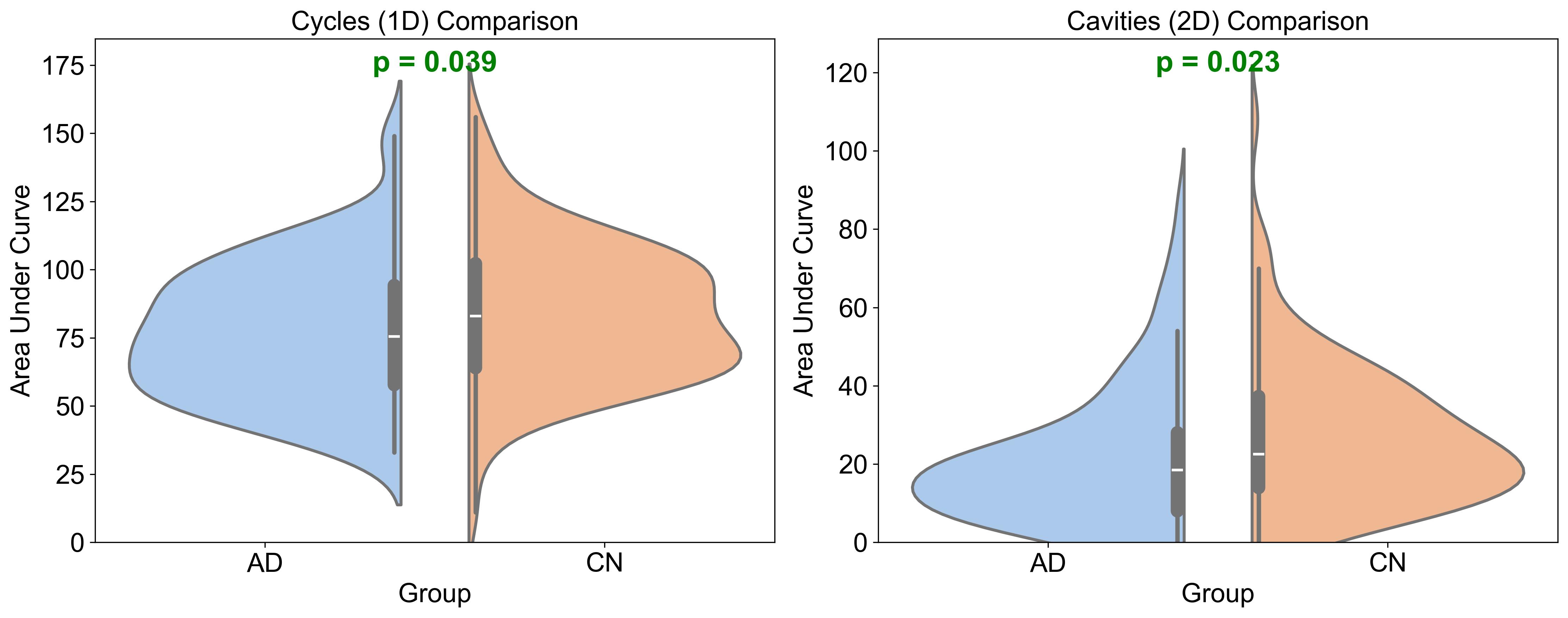}
    \caption{Comparison of higher-order topological features between AD and CN.}
    \label{fig5}
\end{figure}

We computed the area under the Betti curves obtained in Section \ref{quantification of Persistent Topological Features} to quantify the average number of cycles and cavities that emerge during the persistent homology process. Analysis of brain networks from patients with Alzheimer’s disease and healthy controls revealed significant group differences in these features, with all comparisons reaching high statistical significance ($p < 0.05$). Specifically, healthy controls exhibited a greater number of cycles and cavities (Figure \ref{fig5}). These findings are consistent with existing knowledge of Alzheimer’s disease \citep{49}. 

From a neuroscience perspective, the higher number of cycles in healthy individuals is critical for information processing, as cycles reflect feedback loops in the brain that support cyclic transmission, reinforcement, and refinement of information. Such loops help maintain stability and adaptability, enabling efficient processing of complex information. In contrast, the reduction of cycles in Alzheimer’s patients indicates impaired feedback mechanisms, compromising the brain’s ability to integrate and regulate information. For instance, during cognitive control, the integrity of cyclic structures is essential for sustaining attention and suppressing distractions. A reduction in cycle numbers may therefore contribute to symptoms such as attention deficits and reduced cognitive flexibility in Alzheimer’s patients.

Although cavities in brain networks are more abstract, they also carry essential information. The higher number of cavities in healthy controls reflects a well-preserved hierarchical organization and functional segregation within brain networks. Cavities facilitate spatial organization, enabling the isolation and coordination of distinct functional modules. In Alzheimer’s patients, the observed reduction in cavity numbers suggests disruption of this hierarchical organization, impairing inter-regional collaboration and ultimately contributing to cognitive decline.

Overall, differences in higher-order persistent topological features not only provide potential biomarkers for the early diagnosis of Alzheimer’s disease but also offer important insights into the pathological mechanisms underlying the disorder.

\subsubsection{Interpretability of Output Embeddings}
\begin{figure}[t]
\centering

\subfigure{
\includegraphics[scale=0.25]{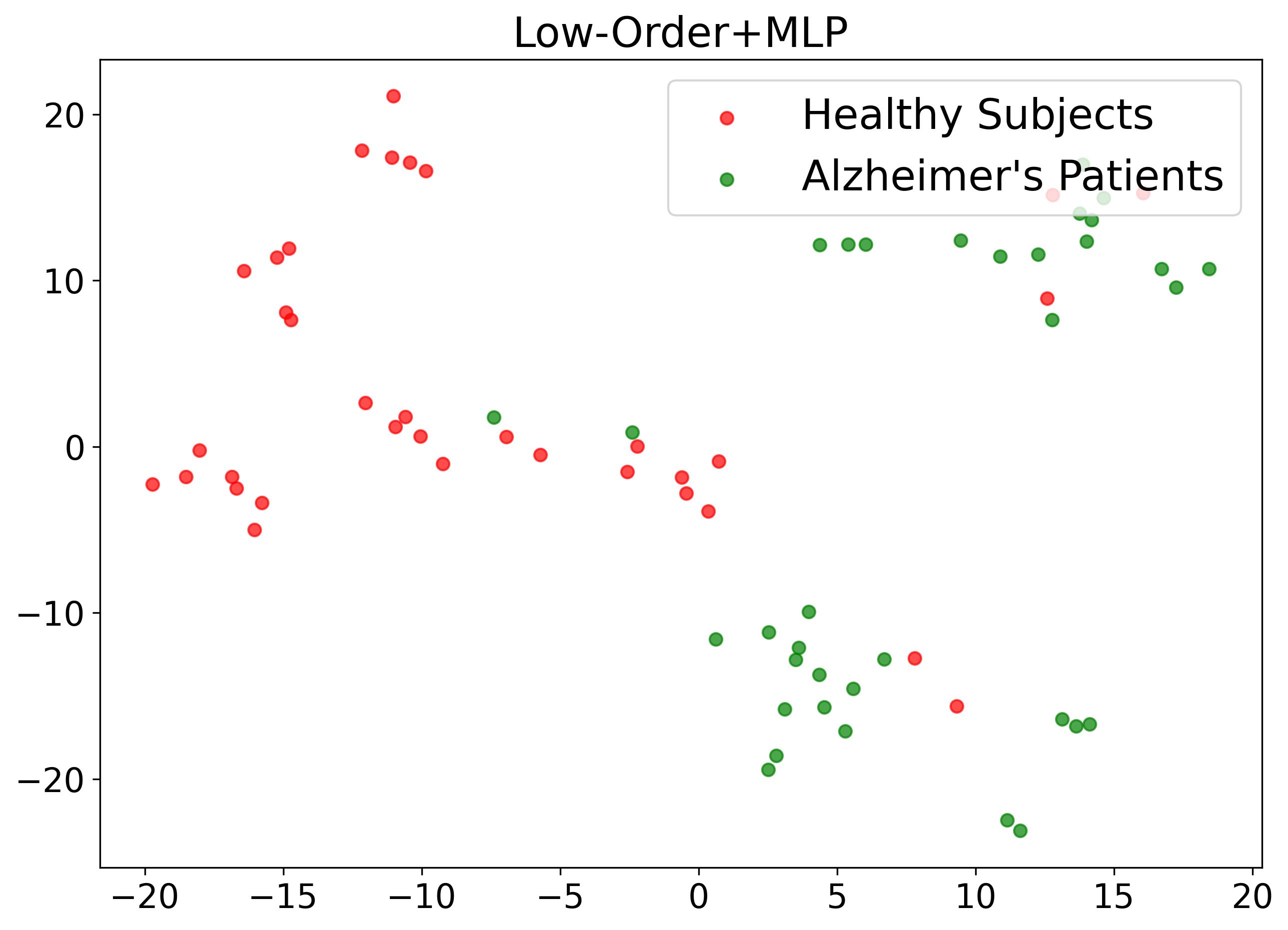}
}
\subfigure{
\includegraphics[scale=0.25]{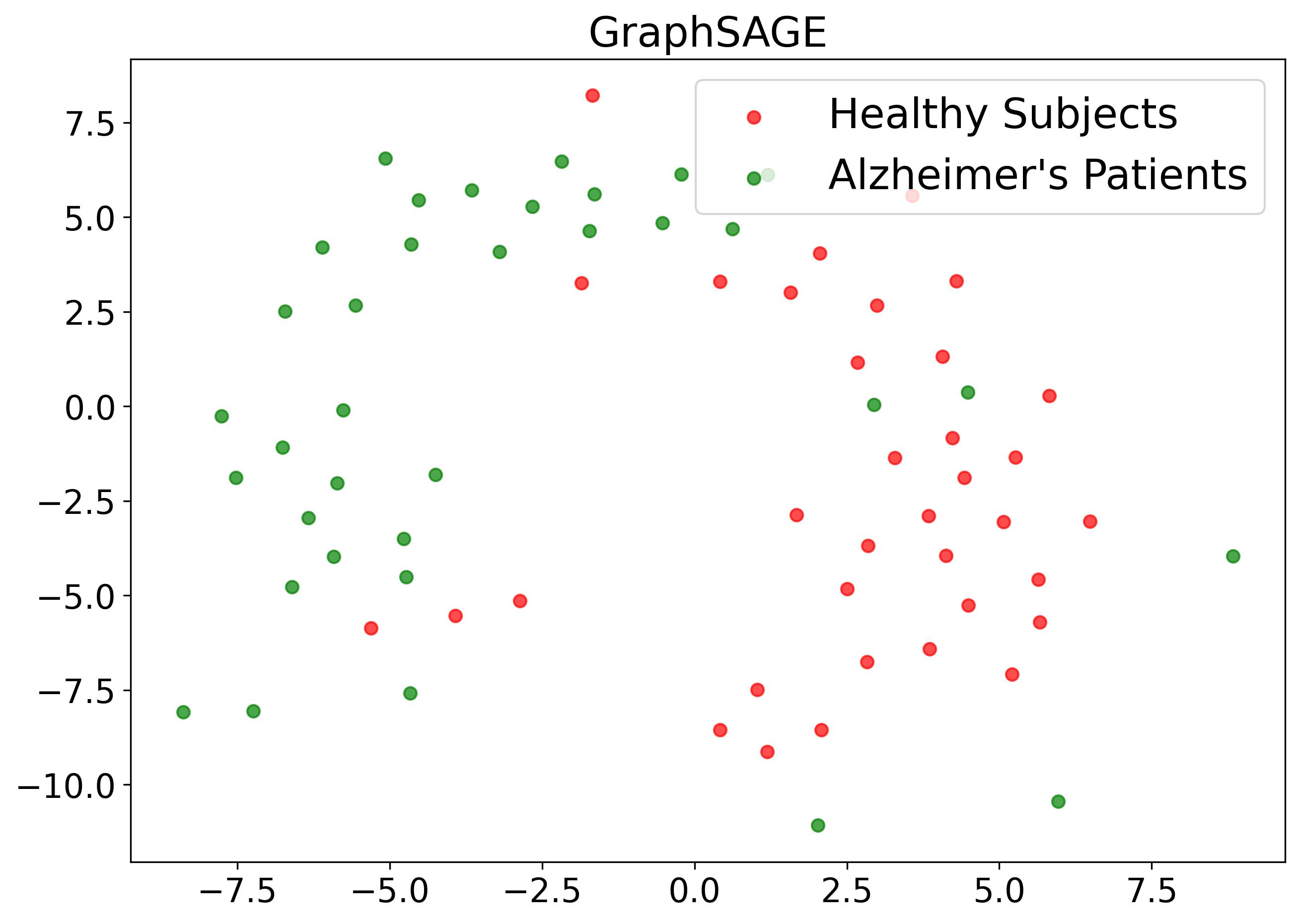}
}

\subfigure{
\includegraphics[scale=0.25]{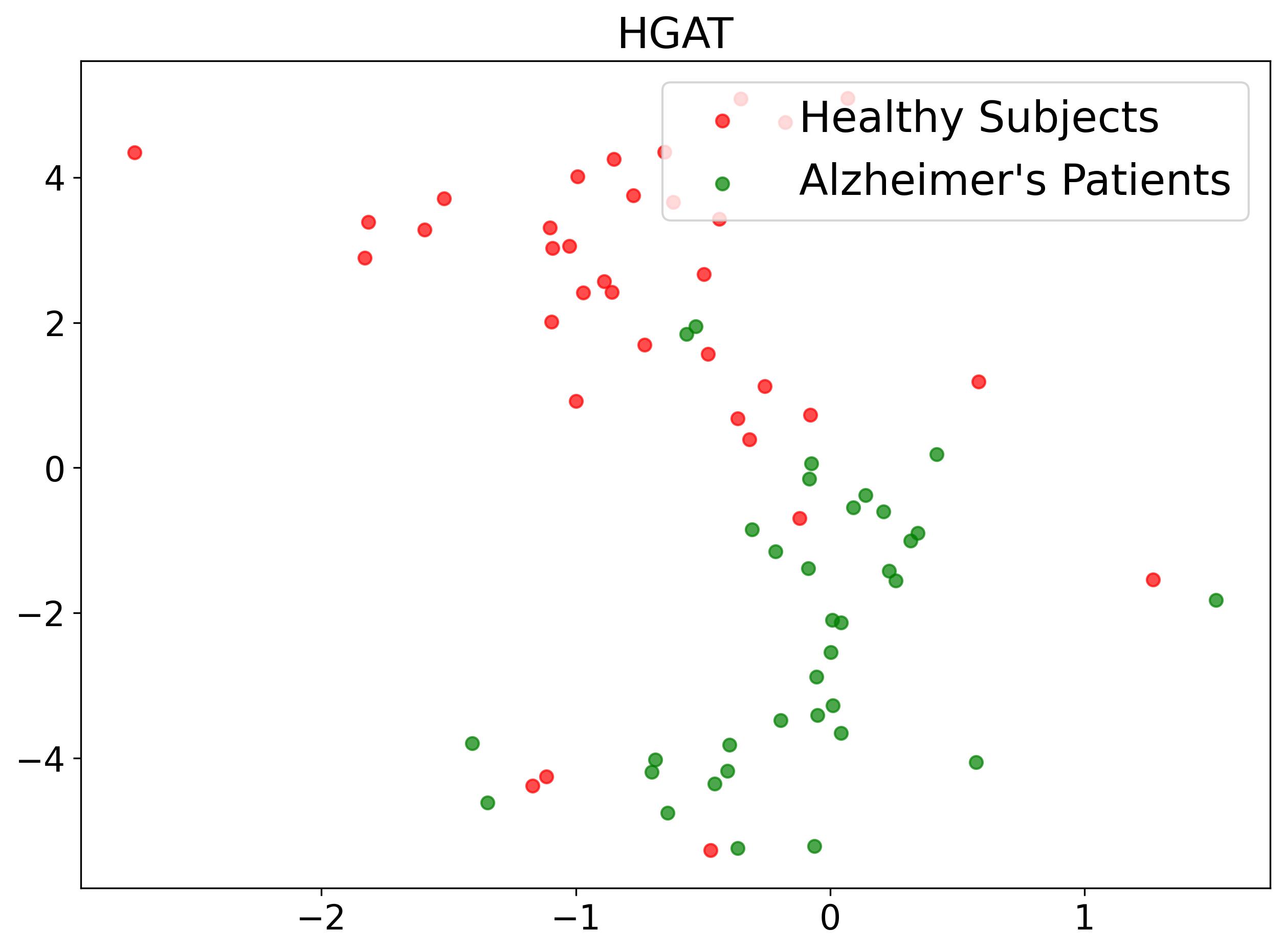}
}
\subfigure{
\includegraphics[scale=0.25]{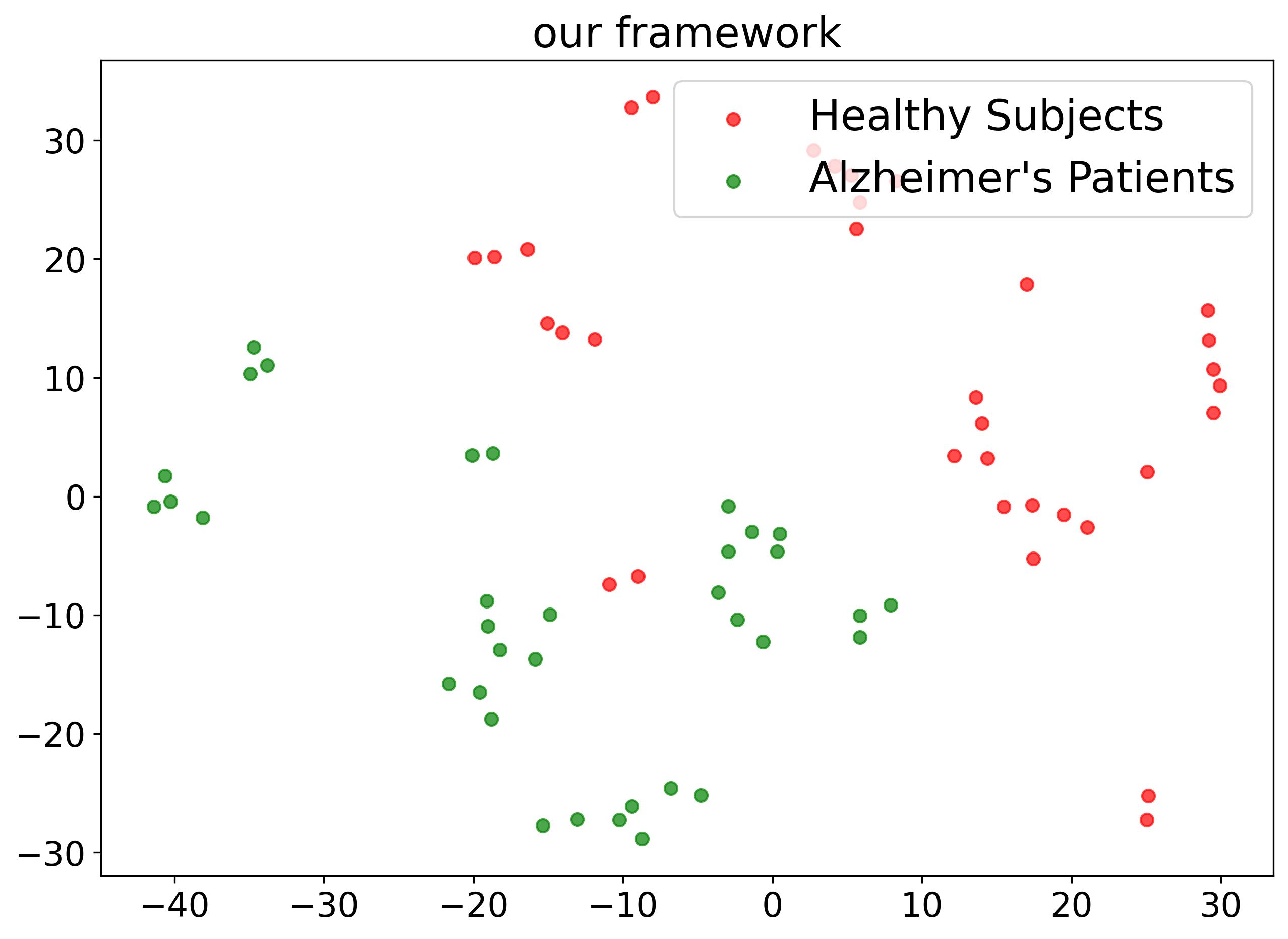}
}
\caption{The t-SNE visualization of the output embeddings learned by Lower-order+MLP, GraphSAGE, HGAT and our framework, respectively.}
\label{fig6}
\end{figure}

To evaluate the interpretability of the proposed characterization method, we visualized the test data from the ADNI1 dataset. Specifically, we employed the Higher-order+MLP model and extracted the output embeddings from the final layer prior to the SoftMax function. These embeddings were visualized using t-SNE \citep{50}, as shown in Figure \ref{fig6}. In the figure, red dots denote healthy subjects, while green dots represent patients with Alzheimer's disease. As illustrated in Figure \ref{fig6}, the results obtained with Lower-order+MLP, GraphSAGE, and HGAT are unsatisfactory, as nodes with different labels are intermixed and lack clear separation. In contrast, the proposed method demonstrates superior performance, with the learned embedding features forming distinct and well-defined boundaries between patients and healthy controls. This experiment confirms that the higher-order persistent topological feature characterization method enables clear discrimination between groups, thereby exhibiting strong interpretability.

\subsubsection{Interpretable Identification of Key Brain Regions}
To identify brain regions significantly associated with Alzheimer's disease, we designed an experiment emphasizing higher-order topological structures, particularly cycles and cavities in brain networks. Using Betti curves with 100 thresholds, we computed t-statistics to quantify differences between AD patients and CN.  
\begin{table}
\centering
\caption{Top ten thresholds showing the most significant differences in cycle and cavity frequencies between brain networks of Alzheimer's disease patients and healthy controls.}
\label{tab:top_thresholds}
\scalebox{0.85}{
\begin{tabular}{cc}
\toprule
Threshold &  t-statistic\\
\midrule
1.30 & 7.768148 \\
1.24 & 5.569526 \\
1.28 & 4.555588 \\
0.78 & 4.378891 \\
1.26 & 4.152438 \\
0.88 & 3.681495 \\
1.34 & 3.663030 \\
1.32 & 3.644254 \\
0.76 & 2.652714 \\
0.70 & 2.260419 \\
\bottomrule
\end{tabular}
}
\end{table}

The ranked t-statistics identified the ten thresholds at which the AD and CN groups exhibited the most significant differences in cycle and cavity counts (see Table \ref{tab:top_thresholds}). These thresholds highlight the networks most affected by AD-related alterations in topological features.

From these, we selected five representative thresholds (0.70, 0.78, 0.88, 1.24, 1.30). For each threshold, we constructed a weighted simplicial complex where edges were included if their weights were smaller than the specified threshold, and then evaluated node importance based on their participation in cycles and cavities. To quantify this, we computed the cocycles of each weighted simplicial complex, thereby capturing the boundary information of \textit{k}-dimensional holes. This boundary information allowed us to derive the number of higher-order topological features associated with each node.  

\begin{figure}
    \centering
    \includegraphics[scale=0.6]{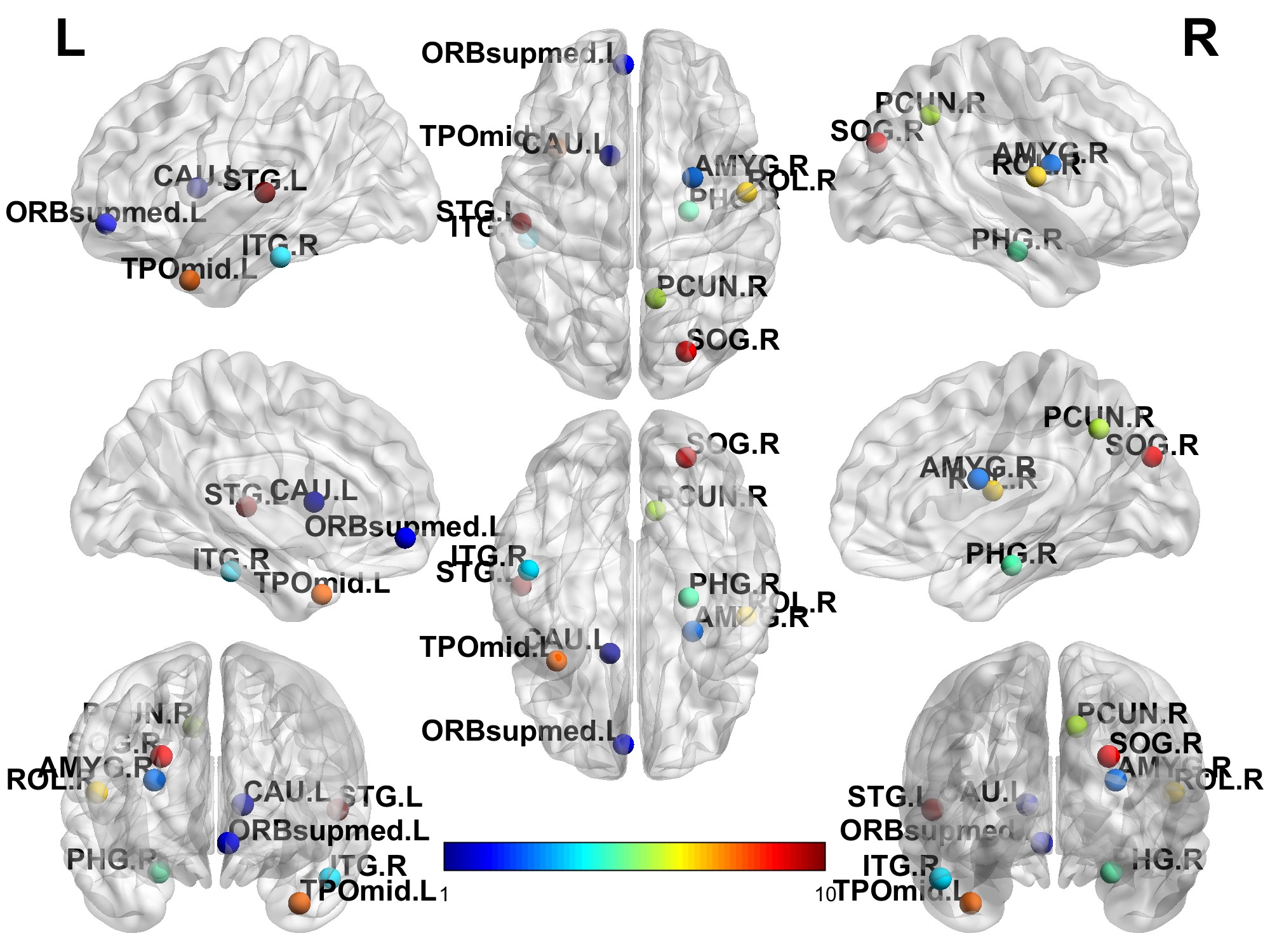}
    \caption{Discriminative brain regions between AD and CN.}
    \label{fig8}
\end{figure}

We then conducted t-tests to identify the top ten nodes with the most significant differences in higher-order topological features between AD and CN groups, with statistical significance set at $p < 0.05$ (Table \ref{tab:top_nodes}). Results were aggregated using a voting mechanism to determine the top ten brain regions consistently exhibiting the most pronounced group differences across thresholds. Table \ref{tab:top_nodes_ad} presents these regions, along with their names and supporting evidence of direct or indirect involvement in AD, obtained from a PubMed literature search. A visual representation of these regions was generated using the BrainNet Viewer toolbox \citep{51}, as shown in Figure \ref{fig8}.

\begin{table}[h]
\centering
\caption{Top ten nodes showing the most significant importance differences between Alzheimer's disease patients and healthy controls across five key filtration thresholds ($\varepsilon_t$).}
\label{tab:top_nodes}
\scalebox{0.85}{
\begin{tabular}{cccccc}
\toprule
\multirow{2}{*}{Rank} & \multicolumn{5}{c}{Node ID at Different Thresholds ($\varepsilon_t$)} \\
\cmidrule(lr){2-6}
 & {$\varepsilon_1=0.7$} & {$\varepsilon_2=0.78$} & {$\varepsilon_3=0.88$} & {$\varepsilon_4=1.24$} & {$\varepsilon_5=1.30$} \\
\midrule
1  & 71 & 71 & 71 & 71 & 71 \\
2  & 25 & 25 & 25 & 25 & 25 \\
3  & 42 & 42 & 42 & 42 & 42 \\
4  & 90 & 90 & 90 & 90 & 90 \\
5  & 18 & 40 & 40 & 40 & 40 \\
6  & 60 & 18 & 68 & 68 & 68 \\
7  & 72 & 68 & 18 & 18 & 18 \\
8  & 68 & 81 & 87 & 87 & 87 \\
9  & 35 & 87 & 50 & 50 & 50 \\
10 & 81 & 50 & 81 & 81 & 81 \\
\bottomrule
\end{tabular}
}
\end{table}

\begin{table}[ht]
\centering
\caption{Top ten nodes with the most significant differences in importance between Alzheimer's disease patients and healthy controls.}
\label{tab:top_nodes_ad}
\scalebox{0.85}{
\begin{tabular}{cll}
\toprule

\textbf{Node} & \textbf{Brain Region} & \textbf{AD Association} \\
\midrule
71 & Caudate (CAU.L) & Direct \cite{52} \\
25 & Superior Medial Orbital Gyrus (ORBsupmed.L) & Indirect \citep{58} \\
42 & Amygdala (AMYG.R) & Direct \citep{57} \\
90 & Inferior Temporal Gyrus (ITG.R) & Direct \citep{56} \\
40 & Parahippocampal Gyrus (PHG.R) & Direct \citep{55} \\
68 & Precuneus (PCUN.R) & Direct \citep{54} \\
18 & Rolandic Operculum (ROL.R) & Direct \citep{52} \\
87 & Middle Temporal Pole (TPOmid.L) & Direct \citep{53} \\
50 & Superior Occipital Gyrus (SOG.R) & Indirect \citep{58} \\
81 & Superior Temporal Gyrus (STG.L) & Indirect \citep{58} \\
\bottomrule
\end{tabular}
}
\end{table}

Our method identified several critical brain regions significantly associated with Alzheimer's disease pathology, encompassing areas involved in memory integration, language processing, cognitive control, and visual recognition. These findings are consistent with established neuroscientific knowledge and supported by extensive literature. The Rolandic Operculum (ROL.R), which plays a key role in speech and language processing, has been shown to accumulate substantial $\beta$-amyloid deposits in AD patients \citep{52}. The left middle temporal pole (TPOmid.L), essential for memory integration and language processing, exhibits altered connectivity patterns in differential methylation region (mDMR) analyzes \citep{53}. The amygdala–prefrontal circuitry (AMYG.R) displays both functional and microstructural abnormalities that contribute to AD pathology, serving as a neural substrate for emotional dysregulation in patients \citep{57}. The parahippocampal gyrus (PHG.R), critical for spatial memory, demonstrates multiple AD-related impairments, including time-dependent inaccuracies in memory-guided eye movements and deficits in long-term spatial memory \citep{55}. In the precuneus (PCUN.R), novelty-related activity follows an inverted U-shaped trajectory across disease progression: elevated in subjective cognitive decline (SCD) and mild cognitive impairment (MCI) compared to CN, but reduced in AD patients relative to MCI \citep{54}. The right inferior temporal gyrus (ITG.R), important for visual processing and semantic memory, undergoes significant synaptic loss during the MCI stage \citep{56}. The left caudate nucleus (CAU.L), implicated in cognitive control and motor coordination, exhibits substantial $\beta$-amyloid deposition, reflecting executive function impairment in AD \citep{52}.

The right superior occipital gyrus (SOG.R), left superior temporal gyrus (STG.L), and left medial orbital gyrus (ORBsupmed.L) are primarily associated with visual–spatial processing, auditory perception, and reward processing, respectively. While direct evidence of AD pathology in these regions is limited, their involvement may reflect secondary network effects or disruptions in functional connectivity \citep{58}. Further investigation of these regions could provide insights into how network-level alterations contribute to disease progression and the heterogeneity of clinical symptoms in AD.

\section{Conclusions}

In this study, we proposed a higher-order topological feature extraction framework capable of capturing connected components, cycles, and cavities across multiple scales in brain networks based on persistent homology. This approach effectively reveals higher-order structures with neurobiological relevance, even in complex and noisy fMRI data. By moving beyond traditional graph-based methods, our framework leverages nested simplicial complex networks to characterize higher-order topological features, thereby overcoming the limitations of fixed-scale brain network analyzes that often discard valuable structural information.
To address the challenges of quantifying high-dimensional features with strong nonlinearity and randomness, we integrated four complementary techniques: Persistence landscape, Betti curve, Heat kernel, and Persistent entropy. These methods enable the effective quantification of complex topological features while enhancing the interpretability of brain connectivity patterns.
Extensive experiments on the ADNI dataset validated the effectiveness of our framework through performance comparisons and ablation studies. Hyperparameter analysis and visualization further confirmed its robustness and stability. Moreover, the identification of key brain regions derived
from cycles and cavities  improved the interpretability of the results, offering insights into disease-related alterations in brain topology. 

Overall, these findings demonstrate the potential of higher-order topological feature extraction in advancing early Alzheimer's disease detection and highlight its broader applicability in the study of neurodegenerative disorders.

\section{Declaration of Competing Interest}
The authors declare that they have no known competing financial interests or personal relationships that could have appeared to
influence the work reported in this paper.
\section{Acknowledgements}
This work was supported in part by the National Natural Science Foundation of China Grant No.12471330 and the National Natural Science Foundation of China Grant No.12231018.
%% For citations use: 
%%       \citet{<label>} ==> Lamport (1994)
%%       \citep{<label>} ==> (Lamport, 1994)
%%

%% The Appendices part is started with the command \appendix;
%% appendix sections are then done as normal sections

%% If you have bib database file and want bibtex to generate the
%% bibitems, please use
%%
\bibliographystyle{elsarticle-harv.bst} 
\bibliography{Reference.bib}

%% else use the following coding to input the bibitems directly in the
%% TeX file.

%% Refer following link for more details about bibliography and citations.
%% https://en.wikibooks.org/wiki/LaTeX/Bibliography_Management

%\begin{thebibliography}{00}

%% For authoryear reference style
%% \bibitem[Author(year)]{label}
%% Text of bibliographic item

\end{document}